\documentclass[11pt]{amsart}

\usepackage{amssymb,amsmath}
\usepackage{bbm}
\usepackage{graphics}
\usepackage{a4wide}

\newtheorem{theorem}{Theorem}
\newtheorem{prop}{Proposition}
\newtheorem{lemma}{Lemma}
\newtheorem{coro}{Corollary}
\newtheorem{fact}{Fact}
 
{
\theoremstyle{definition}

\newtheorem{remark}{Remark}

}

\newcommand{\ts}{\hspace{0.5pt}}

\newcommand{\RR}{\mathbb{R}\ts}
\newcommand{\QQ}{\mathbb{Q}\ts\ts}
\newcommand{\ZZ}{\mathbb{Z}}
\newcommand{\NN}{\mathbb{N}}
\newcommand{\HH}{\mathbb{H}\ts}
\newcommand{\II}{\mathbb{I}}
\newcommand{\oo}{{\displaystyle\mathfrak{o}}}

\newcommand{\alat}{L}
\newcommand{\up}{\!\!\widetilde{\hphantom{aa}}}
\newcommand{\uu}{{\ts\ts\!\times\!}}

\newcommand{\ii}{\mathrm{i}}
\newcommand{\jj}{\mathrm{j}}
\newcommand{\kk}{\mathrm{k}}

\newcommand{\gG}{\varGamma}
\newcommand{\gL}{\varLambda}

\DeclareMathOperator{\lcm}{lcm}
\DeclareMathOperator{\den}{den}
\DeclareMathOperator{\tr}{tr}
\DeclareMathOperator{\nr}{nr}
\DeclareMathOperator{\N}{N}

\DeclareMathOperator{\OC}{\mathrm{OC}}
\DeclareMathOperator{\SOC}{\mathrm{SOC}}
\DeclareMathOperator{\OS}{\mathrm{OS}}
\DeclareMathOperator{\SOS}{\mathrm{SOS}}

\newcommand{\lbar}[1]{\,\overline{\! #1 \!}\,}

\begin{document}

\title[Coincidence rotations of $A_4$]
{Coincidence rotations of the root lattice $A_4$}

\author{Michael Baake}
\author{Uwe Grimm}
\author{Manuela Heuer}
\author{Peter Zeiner}

\address{Fakult\"at f\"ur Mathematik, Universit\"at Bielefeld, \newline
\hspace*{12pt}Postfach 100131, 33501 Bielefeld, Germany}
\email{$\{$mbaake,pzeiner$\}$@math.uni-bielefeld.de}

\address{Department of Mathematics, The Open University, \newline
\hspace*{12pt}Walton Hall, Milton Keynes MK7 6AA, UK}
\email{$\{$u.g.grimm,m.heuer$\}$@open.ac.uk}

\begin{abstract} 
  The coincidence site lattices of the root lattice $A_4$ are
  considered, and the statistics of the corresponding coincidence
  rotations according to their indices is expressed in terms of a
  Dirichlet series generating function.  This is possible via an
  embedding of $A_4$ into the icosian ring with its rich arithmetic
  structure, which recently \cite{BHM} led to the classification
  of the similar sublattices of $A_4$.
\end{abstract}

\maketitle

\bigskip
\centerline{\large Dedicated to Ludwig Danzer on the occasion of his
80th birthday}
\bigskip

\section{Introduction and general setting}

Consider a lattice $\gG$ in Euclidean $d$-space, i.e., a cocompact
discrete subgroup of $\RR^d$. An element $R\in\mathrm{O}(d,\RR)$ is
called a (linear) \emph{coincidence isometry} of $\gG$ when $\gG$ and
$R\gG$ are commensurate, written as $\gG\sim R\gG$, which means that
they share a common sublattice. The intersection $\gG\cap R\gG$ is
then called the \emph{coincidence site lattice} (CSL) for the isometry
$R$.  When this is the case, the corresponding \emph{coincidence
  index} $\varSigma (R)$ is defined as
\[
    \varSigma (R) \, = \, [ \gG : (\gG\cap R\gG)]\, ,
\]
and it is set to $\infty$ otherwise.
The index satisfies $[ \gG : (\gG\cap R\gG)]=[R \gG : 
(\gG\cap R\gG)]$, as $\gG$ and $R\gG$ possess fundamental domains
of the same volume, compare \cite{Cassels} for general background
on lattice theory. Moreover,
\[
    \OC (\gG) \, := \,
    \{ R\in\mathrm{O}(d,\RR)\mid \varSigma (R) < \infty \}
\]
is a group, and one also has $\varSigma (R^{-1})=\varSigma (R)$, see
\cite{B} for a general survey and several typical examples. The
subgroup $\SOC (\gG)$ consists of all rotations within
$\OC (\gG)$.

Coincidence site lattices play an important role in crystallography,
in the description and understanding of grain boundaries, compare
\cite{B} and references given there. In recent years, they have also
found applications in lattice discretisation problems \cite{BLS}. From
a more theoretical angle, they show up in various lattice and tiling
problems, such as Danzer's `Ein-Stein-Tiling' \cite{D1,D2,BF} or the
analysis of the pinwheel tilings of the plane \cite{MPS,BFG}. Apart
from that, several attempts have been made to get further insight into
the theory, see \cite{BG,BPR,Z1,Zou1} and references given there for
recent publications.

Another relevant object in this context, with
$\RR_{+}:=\{\alpha\in\RR\mid \alpha >0\}$, is the set
\[
    \OS (\gG) \, := \, 
    \{ R\in\mathrm{O}(d,\RR)\mid \alpha R\gG\subset\gG
    \text{ for some } \alpha \in\RR_{+} \}\, .
\]
It consists of all linear isometries that emerge from similarity
mappings of $\gG$ into itself, while $\SOS (\gG)$ is then the subset
of rotations of that kind. A sublattice of $\gG$ of the form $\alpha
R\gG$ is called a \emph{similar sublattice} of $\gG$, or SSL for
short.
\begin{fact} \label{os-group}
  If $\gG\subset\RR^d$ is a lattice, $\OS (\gG)$ 
  and\/ $\SOS (\gG)$
  are subgroups of\/ $\mathrm{O} (d,\RR)$.
\end{fact}
\begin{proof}
It suffices to verify the subgroup property for $\OS (\gG)$.  If $R$
and $S$ are isometries in $\OS (\gG)$, with attached positive
numbers $\alpha$ and $\beta$ say, then so is $RS$, because
$\alpha\beta RS\gG = \alpha R (\beta S\gG)\subset\alpha R\gG
\subset\gG$.

Also, one has $\gG\subset\frac{1}{\alpha}R^{-1}\gG$, so that
$[\frac{1}{\alpha}R^{-1}\gG : \gG] = [\gG : \alpha R \gG] =: m\in\NN$.
By standard lattice theory, compare \cite{Cassels}, this means that
$\frac{m}{\alpha}R^{-1}\gG\subset\gG$, which shows that also $R^{-1}$
is an element of $\OS (\gG)$.
\end{proof}

Let us start with a general observation on the connection between the
similarities and the coincidence isometries.
\begin{lemma} \label{lemma-one}
  If $R$ is a coincidence isometry for the lattice $\gG\subset \RR^d$,
  there exists some $\alpha\in\RR_{+}$ so that $\alpha R \gG\subset \gG$.
  In other words, $\OC (\gG)$ is a subgroup of $\ts\OS (\gG)$.
\end{lemma}
\begin{proof}
$R\in\OC (\gG)$ means $\varSigma (R) =
[\gG : (\gG\cap R\gG)] = [R\gG : (\gG\cap R\gG)] = n \in \NN$. This
implies $n R\gG \subset (\gG\cap R\gG) \subset \gG$ by standard
lattice theory.
\end{proof}

The converse is not true in general, meaning that \emph{not} all
rotations and reflections from similar sublattices will give rise to
coincidence isometries. It is precisely one of our goals later on to
find the distinction for the case of the root lattice $A_4$.

In view of Lemma~\ref{lemma-one}, it is reasonable to define the
\emph{denominator} of a matrix $R\in\OS (\gG)$ relative to the lattice
$\gG$ as
\begin{equation}  \label{def-denom}
   \den^{}_{\gG} (R) \, = \,
   \min\ts \{\alpha \in  \RR_{+} \mid \alpha R\gG\subset\gG\}\, .
\end{equation}
Clearly, as $R$ is an isometry, one always has $\den^{}_{\gG} (R)\ge 1$,
and from $\den^{}_{\gG} (R)\ts R\gG\subset\gG$ one concludes that
$\bigl(\den^{}_{\gG} (R)\bigr)^d$ must be an integer.  Consequently,
$\den^{}_{\gG} (R)$ is either a positive integer or an irrational number,
but still algebraic.  Moreover, by standard arguments, one has
\begin{equation} \label{den-values}
  \{ \alpha \in\RR_{+}\mid\alpha R\gG\subset\gG\} \, = \,
  \den^{}_{\gG} (R)\,\NN\, .
\end{equation}
This gives rise to the following refinement of Lemma~\ref{lemma-one}.
\begin{lemma} \label{lemma-two}
  Let $\gG\subset\RR^d$ be a lattice, with groups $\OS (\gG)$ and\/
  $\OC (\gG)$ as defined above.  With the denominator from
  $\eqref{def-denom}$, one has\/ $\OC (\gG) = \{ R\in \OS (\gG)\mid
  \den^{}_{\gG} (R) \in \NN \}$.
\end{lemma}
\begin{proof}
If $\den^{}_{\gG} (R) \in \NN$, one has $\den^{}_{\gG}(R) R\gG\subset
(\gG\cap R\gG)$. Consequently, the lattices $\gG$ and $R\gG$ are
commensurate, so that the inclusion $\{ R\in \OS (\gG)\mid
\den^{}_{\gG} (R) \in \NN \} \subset \OC(\gG)$ is clear.

Conversely, if  $R\in\OC (\gG)$, $\gG$ and $R\gG$  are commensurate by
definition.  In particular, one has $\varSigma (R)\ts R\gG\subset\gG$,
so     that    $\varSigma     (R)\in\den^{}_{\gG}     (R)\ts\NN$    by
Eq.~\eqref{den-values}.   As  $\varSigma   (R)\in\NN$,  this  is  only
possible if  $\den^{}_{\gG} (R)\in\QQ$. Since we also  know from above
that  $\bigl(\den^{}_{\gG} (R)\bigr)^d\in\NN$,  we  may now conclude  that
$\den^{}_{\gG} (R)\in\NN$, whence the claim follows.
\end{proof}

For later use, we state a factorisation property for coincidence
indices from \cite{BZ2}.
\begin{lemma}   \label{factor-lemma}
  Let $\gG\subset\RR^d$ be a lattice and $R_1 , \ts R_2 \in \OC
  (\gG)$.  When $\varSigma (R_1)$ and $\varSigma (R_2)$ are relatively
  prime, one has $\varSigma (R_1\ts R_2) = \varSigma (R_1)\,\varSigma
  (R_2)$.  In general, one has the divisibility relation $ \varSigma
  (R_1\ts R_2) \, | \, \varSigma (R_1)\,\varSigma (R_2)$.  \qed
\end{lemma}

When a lattice $\gG$ is given, the set $\varSigma (\OC (\gG))$ is
called the \emph{simple coincidence spectrum}. It may or may not
possess an algebraic structure. In nice situations, $\varSigma (\OC
(\gG))$ is a multiplicative monoid within $\NN$.  On top of the
spectrum, one is also interested in the number $f(m)$ of different
CSLs of a given index $m$. This arithmetic function is often
encapsulated into a Dirichlet series generating function,
\begin{equation} \label{def-diri}
   \varPhi_{\gG} (s) \, := \, 
    \sum_{m=1}^{\infty} \frac{f(m)}{m^s} \, ,
\end{equation}
which is a natural approach because it permits an Euler product
decomposition when $f$ is multiplicative. Beyond three dimensions,
this generating function seems difficult to determine, and $A_4$ is no
exception. In this paper, we thus concentrate on the slightly simpler
problem to count the coincidence rotations of $A_4$ of index $m$,
which always come in multiples of $120$, the order of the rotation
symmetry group of $A_4$. The latter is the subgroup of $\SOC (A_4)$ of
rotations with coincidence index $\varSigma=1$. If $120 f_{\rm rot}
(m)$ is the number of coincidence rotations with index $m$, we define
the corresponding generating function as
\begin{equation} \label{def-diri-rot}
   \varPhi_{\gG}^{\rm \, rot} (s) \, := \, 
    \sum_{m=1}^{\infty} \frac{f_{\rm rot}(m)}{m^s} \, ,
\end{equation}
which will turn out to possess a nice Euler product expansion.  Note
that $0 \le f(m) \le f_{\rm rot} (m)$ and that $f(m) \neq 0$ if and
only if $f_{\rm rot}(m)\neq 0$.  Asymptotic properties for
$m\to\infty$ can be extracted from the generating functions via
Delange's theorem, see \cite{BM} and references therein for details in
this context.

\section{The root lattice $A_4$ and its arithmetic structure}

\begin{figure}
 \scalebox{0.6}{\includegraphics{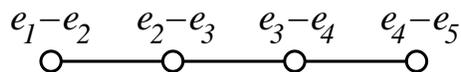}}
 \caption{Standard basis representation of the root lattice $A_4$.
 \label{a4-fig}}
 \end{figure}

The root lattice $A_4$ is usually defined as a lattice in a
$4$-dimensional hyperplane of $\RR^5$, via the Dynkin diagram of
Figure~\ref{a4-fig}. Here, the $e_i$ denote the standard Euclidean
basis vectors in $5$-space. Though convenient for many purposes, this
description does not seem to be optimal for the geometric properties
we are after.  Since the similar sublattices of $A_4$ were recently
classified \cite{BHM} by a different $4$-dimensional approach, using
the arithmetic of the quaternion algebra $\HH (\QQ(\sqrt{5}\,))$, we
use the same setting again in this paper. 

{}From now on, we use the notation $K=\QQ(\sqrt{5}\,)$ for brevity.
The algebra $\HH (K)$, which is a skew field, is explicitly given as
$\HH (K) = K \oplus \ii K \oplus \jj K \oplus \kk K$, where the
generating elements satisfy Hamilton's relations $\ii^2 = \jj^2 =
\kk^2 = \ii\jj\kk = -1$, see \cite{O} for more. $\HH (K)$ is equipped
with a \emph{conjugation} $\,\bar{.}\,$ which is the unique mapping that
fixes the elements of the centre of the algebra $K$ and reverses the
sign on its complement. If we write $q=(a,b,c,d)=a+\ii b+\jj c+\kk d$,
this means $\bar{q}=(a,-b,-c,-d)$.

The reduced norm and trace in $\HH (K)$ are defined as usual
\cite{V,MW,O} by
\begin{equation} \label{norm-trace}
   \nr (q) \, = \, q\bar{q} \, = \, \lvert q\rvert^2 
   \quad \text{and} \quad
   \tr (q) \, = \, q + \bar{q} \, ,
\end{equation}
where we canonically identify an element $\alpha\in K$ with the
quaternion $(\alpha,0,0,0)$.  For any $q\in\HH (K)$, $\lvert q \rvert$ is
its Euclidean length, which need not be an element of $K$.
Nevertheless, one has $\lvert rs\rvert = \lvert r\rvert \lvert
s\rvert$ for arbitrary $r,s\in\HH (K)$.  Due to the geometric meaning,
we use the notations $\lvert q\rvert^2$ and $\nr (q)$ in parallel.  An
element $q\in\HH (K)$ is called \emph{integral} when both $\nr (q)$
and $\tr (q)$ are elements of $\oo:=\ZZ[\tau]$, which is the ring of
integers of the quadratic field $K$, where $\tau = (1+\sqrt{5}\,)/2$
is the golden ratio.

In this setting, we use the lattice
\begin{equation} \label{def-alat}
   \alat \, = \, 
     \bigl\langle (1,0,0,0), \tfrac{1}{2}(-1,1,1,1), (0,-1,0,0), 
       \tfrac{1}{2}(0,1,\tau\!-\!1,-\tau) \bigr\rangle_{\ZZ}\, ,
\end{equation}
which is the root lattice $A_4$ relative to the inner product $\tr
(x\bar{y}) = 2 \langle x\ts | \ts y\rangle$, where
$\langle x\ts | \ts y\rangle$ denotes the standard Euclidean
inner product, see \cite{CMP,BHM} for
details. This way, $\alat$ is located within the icosian ring $\II$,
\begin{equation} \label{def-icosian}
   \II \, = \, \bigl\langle (1,0,0,0),(0,1,0,0),
   \tfrac{1}{2}(1,1,1,1),\tfrac{1}{2}(1 \! - \!\tau,\tau,0,1)
   \bigr\rangle_{\oo}\, ,
\end{equation}
which is a maximal order in $\HH (K)$, compare \cite{MP,MW} and
references given there. In particular, all elements of $\II$ (and
hence also those of $\alat$) are integral in $\HH (K)$.  In fact, one
can use the quadratic form defined by $\tr (x\bar{y})$ to define the
\emph{dual} of a full $\oo$-module $\gL\subset \HH (K)$ as
\begin{equation} \label{dual-module}
    \gL^* \, = \, \{ x\in\HH (K) \mid \tr (x\bar{y})\in\oo
     \text{ for all } y\in\gL \} \, .
\end{equation}
With this definition, one has the following important property
of the icosian ring, compare \cite{S,MP,CMP} for details.
\begin{fact} \label{self-dual}
   The icosian ring is self-dual, i.e., one has $\II^* = \II$.
   \qed
\end{fact}

Since $\HH (K)$ has class number $1$, compare \cite{S,V}, all ideals
of $\II$ are principal. The detailed arithmetic structure of $\II$
was the key to solving the related sublattice problem \cite{BHM} for
$\alat$. What is more, one significantly profits from another
map, called \emph{twist map} in \cite{BHM}, which is an involution of
the second kind for $\HH (K)$. If $q=(a,b,c,d)$, it is defined by the
mapping $q\mapsto \widetilde{q}$,
\begin{equation} \label{def-twist}
    \widetilde{q} \, = \, (a{\ts}',b{\ts}',d{\ts}',c{\ts}')\, ,
\end{equation}
where ${}'$ denotes algebraic conjugation in $K$, as defined by the
mapping $\sqrt{5}\mapsto -\sqrt{5}$. Note the unusual combination
of algebraic conjugation of all coordinates with a permutation of
the last two -- which also explains the choice of the term
`twist map'. The algebraic conjugation in $K$ is
also needed to define the absolute norm on $K$, via 
\[
     \N (\alpha) \, = \, \lvert \alpha \alpha{\ts}' \rvert\,  .
\]

{}For the various properties of the twist map, we refer the reader to
\cite{BHM} and references therein.  The most important one in our
present context is the relation between $\alat$ and $\II$.
\begin{fact} \label{fundamental}
   Within $\HH (K)$, one has\/ $\widetilde{\II} = \II$ and 
   $\alat = \{ x\in\II\mid\widetilde{x}=x\}$. \qed
\end{fact}

Another useful characterization is possible via the $\QQ$-linear
mappings $\phi^{}_{\pm} \! : \, \HH (K) \longrightarrow
\HH (K)$, defined by $\phi^{}_{\pm} (x) = x\pm\widetilde{x}$,
which are connected via the relation
\[
     \phi^{}_{\pm} (\sqrt{5}\ts x) \, = \, 
     \sqrt{5}\,\phi^{}_{\mp} (x)
\]
and the obvious property $\mathrm{ker} (\phi^{}_{\pm}) =
\mathrm{im} (\phi^{}_{\mp})$.
\begin{lemma} \label{new-sum}
  The lattice $\alat\ts$ from $\eqref{def-alat}$ satisfies\/
  $\alat = \{ x + \widetilde{x}\mid x\in\II\} = \phi^{}_{+} (\II)$.
\end{lemma}
\begin{proof}
For any $x\in\II$, we clearly have $x+\widetilde{x}\in\alat$ by 
Fact~\ref{fundamental}. On the other hand, observing
$\tau'=1-\tau$, any $x\in\alat$ permits the decomposition
\[
   x \, = \, (\tau + \tau')\, x \, = \,
   \tau\ts x + \tau' \widetilde{x} \, = \,
   \tau x + \widetilde{\tau x} \, .
\]
Since $x\in\II$ and $\II$ is an $\oo$-module, we still have $\tau
x\in\II$, and the claim follows.
\end{proof}

\begin{remark}
{}For $0\neq q\in\II$, one has $\nr (q)\in\oo$ and
$\nr (q) = q\bar{q} >0$. As also $\widetilde{q}\in\II$,
one finds $\nr (\widetilde{q})\in\oo$ with
$\nr (\widetilde{q})>0$ and $\nr (\widetilde{q}) =
\nr (q)'$, so that $\nr (q)$ is always a totally positive
element of $\oo$.
\end{remark}

\section{Coincidence site lattices via quaternions}

It is clear that we can restrict ourselves to the investigation of
rotations only, because $\overline{\alat}=\alat$, so any orientation
reversing operation can be obtained from an orientation preserving one
after applying conjugation first.

Let us start by recalling a fundamental result from \cite{BHM}.
\begin{fact} \label{prop-one}
  If $q\in\II$, one has $q\alat\widetilde{q}\subset\alat$.  Moreover,
  all similar sublattices of $\alat$ are of the form $q \alat
  \widetilde{q}$ with $q\in\II$. \qed
\end{fact}

{}For a given SSL of $\alat$, now written as $q\alat\widetilde{q}$,
the corresponding rotation is then given by the mapping $x\mapsto
\frac{1}{\lvert q \widetilde{q}\ts\rvert}\,q x \widetilde{q}$.  It is
clear that many different $q$ result in the same rotation.  In fact,
we can restrict $q$ to suitable subsets of icosians without missing
any rotation, which we shall do later on.

Below, we need a refinement of Fact~\ref{prop-one}. Recall that
a sublattice $\gL$ of $\alat$ is called \emph{$\alat$-primitive} when
$\alpha\gL\subset\alat$, with $\alpha\in\QQ$, implies $\alpha\in\NN$.
Similarly, an element $p\in\II$ is called \emph{$\II$-primitive} when
$\alpha p\in\II$, this time with $\alpha\in K$, is only possible with
$\alpha\in\oo$, see \cite{BHM} for details, and for a proof of
the following result. For brevity, we simply use the term ``primitive''
in both cases, as the meaning is clear from the context.
\begin{prop} \label{prop-one-prime}
  A similar sublattice of $\alat$ is primitive if and only if
  it is of the form $q\alat\widetilde{q}$ with $q$ a primitive
  element of $\ts\II$. \qed
\end{prop}

By Lemmas~\ref{lemma-one} and \ref{lemma-two}, we know how 
$\SOC (\alat)$ and $\SOS (\alat)$ are related in general. 
Here, Fact~\ref{prop-one} tells us that any
similarity rotation of $\alat$ is of the form $x\mapsto
\frac{1}{\lvert q\widetilde{q}\ts\rvert}\, q x \widetilde{q}$ with
$q\in\II$. Among these, we have to identify the $\SOC (\alat)$
elements, which is possible as follows.
\begin{coro} \label{commensurate}
  Let $0\neq q\in\II$ be an arbitrary icosian.  The lattice
  $\frac{1}{\lvert q\widetilde{q}\ts\rvert}\, q \alat \widetilde{q}$
  is commensurate with $\alat$ if and only if\/ $\lvert q
  \widetilde{q}\ts\rvert\in\NN$.
\end{coro}
\begin{proof}
When $q=\alpha r$ with $0\neq \alpha\in \oo$, one has $\lvert
q\tilde{q}\rvert = \N (\alpha) \lvert r\tilde{r}\rvert$ with $\N
(\alpha)\in \NN$. If $q$ is primitive, the claim is clear by
Lemma~\ref{lemma-two} because $\lvert q\tilde{q}\rvert$ is then the
denominator of the rotation $x\mapsto \frac{1}{\lvert
q\widetilde{q}\ts\rvert}\, q x \widetilde{q}$. Otherwise, $q$ is an
$\oo$-multiple of a primitive icosian, $r$ say, and the claim follows 
from the initial remark.
\end{proof}

Let us call an icosian $q\in\II$ \emph{admissible} when
$\lvert q \widetilde{q}\ts\rvert\in\NN$. As $\nr (\widetilde{q}\ts) =
\nr (q){\ts}'$, the admissibility of $q$ implies that
$\N\bigl(\nr (q)\bigr)$ is a square in $\NN$.

\begin{theorem}
  The CSLs of $\alat$ are precisely the lattices of the form $\alat
    \cap \frac{1}{\lvert q\widetilde{q}\ts\rvert}\, q\alat\widetilde{q}$
    with $q\in\II$ primitive and admissible.
\end{theorem}
\begin{proof}
All CSLs can be obtained from a rotation, as $\alat$ is invariant
under the conjugation $x\mapsto\bar{x}$. By Lemma~\ref{lemma-one}, we
need only consider rotations from $\OS (\alat)$, which, by
Fact~\ref{prop-one}, are all of the form
$x\mapsto\frac{1}{\lvert y\widetilde{y}\ts\rvert}\,y x\widetilde{y}$
with $y\in\II$. If $y$ is not primitive, we can write it as
$y=\alpha\ts q$ with $\alpha\in\oo$ and $q\in\II$
primitive. Since $\alpha$ is a central element, $y$ and $q$ define the
same rotation. An application of Corollary~\ref{commensurate} now
gives the claim.
\end{proof}

This is the first step to connect certain primitive right
ideals $q\II$ of the icosian ring with the CSLs of $\alat$.
Before we continue in this direction, let us consider the
relation with the coincidence rotations.

\begin{lemma} \label{symmetries} 
  Let $r,s\in\II$ be primitive and admissible quaternions, with
  $r\II=s\II$.  Then, one has $\alat\cap \frac{r\alat\ts
    \tilde{r}}{\lvert r \tilde{r}\rvert} = \alat\cap \frac{s\alat
    \tilde{s}}{\lvert s \tilde{s}\rvert}$.
\end{lemma}
\begin{proof}
  When $r\II = s\II$, one has $s=r\varepsilon$ for some
  $\varepsilon\in\II^{\uu}$, where $\II^{\uu}$ denotes the unit group
  of $\II$, see \cite{MW} for its structure. Since, by
  \cite[Lemma~4]{BHM}, we then know that
  $\varepsilon\alat\ts\widetilde{\varepsilon}= \alat$, one has
  $r\alat\ts \tilde{r} = s\alat \tilde{s}$ in this case. As $\nr
  (\varepsilon\widetilde{\varepsilon}\,) = \N (\nr(\varepsilon))=1$,
  one also finds $\lvert s\tilde{s}\rvert = \lvert r\tilde{r}\rvert$.
  Consequently, $\frac{r\alat\ts \tilde{r}}{\lvert r \tilde{r}\rvert}
  = \frac{s\alat \tilde{s}}{\lvert s \tilde{s}\rvert}$, and the CSLs
  of $\alat$ defined by $r$ and $s$ are equal.
\end{proof}

\begin{remark} \label{weak-converse} The converse statement to
  Lemma~\ref{symmetries} is not true, as the equality of
  two CSLs does \emph{not} imply the corresponding rotations to be
  symmetry related. An example is provided by $r=(\tau,2\tau,0,0)$ and
  $s=(\tau^2,\tau,\tau,1)$, which define the same CSL, though
  $s^{-1}r$ is not a unit in $\II$.  The CSL is spanned by the
  basis $\{(1,2,0,0),(2,-1,0,0),(\frac{3}{2},\frac{1}{2},\frac{1}{2},
  \frac{1}{2}),(-1,\frac{1}{2},\frac{\tau-1}{2},-\frac{\tau}{2})\}$.
  However, when two primitive
  quaternions $r,s\in\II$ define rotations that are related by a
  rotation symmetry of $A_4$, one has $r\II=s\II$ as a direct
  consequence of \cite[Lemma~4]{BHM}.
\end{remark}

Although the primitive elements of $\II$ are important in this
context, we need a variant for our further discussion. Let $q\in\II$
be primitive and admissible. Since $\oo$ is Dedekind, one has the
relation $(x\oo)^{-1} = \frac{1}{x}\oo$ for any principal fractional 
ideal with nonzero $x\in K$, see \cite[Ch.~I.4]{J} for details. Then, 
the fractional ideal
\[
    \bigl(\nr(q)\oo\cap\nr(\widetilde{q}\,)\oo\bigr)^2
    \bigl(\lvert q\widetilde{q}\ts\rvert^2 \oo\bigr)^{-1}
    \, = \, \frac{\bigl(\lcm (\nr (q),\nr (\widetilde{q}\,))\bigr)^2}
     {\lvert q\widetilde{q}\ts\rvert^2}\,\oo 
    \, = \, \beta_{q} \oo\, \beta_{\tilde{q}}\oo
\]
is a square as well, where $\beta_{q} := \lcm (\nr (q),\nr
(\widetilde{q}\,))/\nr (q)\in\oo$ is well-defined up to units of
$\oo$, as $\oo$ is a principal ideal domain. Clearly, $\beta_{q}\oo$
and $\beta_{\tilde{q}}\oo$ are coprime by construction. Since their
product is a square in $\oo$ (up to units), we have $\beta_{q}\oo =
\bigl(\alpha_{q}\oo\bigr)^2$ for some $\alpha_q\in\oo$. Explicitly, we
may choose
\begin{equation} \label{def-alpha}
    \alpha^{}_{q} \, = \, 
    \sqrt{\frac{\lcm(\nr (q),\nr (\widetilde{q}\,))}{\nr (q)}}
    \, = \,
    \sqrt{\frac{\lcm(\nr (q),\nr (q){\ts}')}{\nr (q)}}
    \, \in \, \oo \, ,
\end{equation}
where we assume a suitable standardisation for the $\lcm$ of
two elements of $\oo$. Again, $\alpha^{}_{q}$ is only defined
up to units of $\oo$, which is tantamount to saying that we
implicitly work with the principal ideal $\alpha^{}_{q}\oo$
here. Moreover, we have the relation $\alpha^{}_{\tilde{q}}
= \widetilde{\alpha^{}_{q}} = \alpha^{\,\prime}_{q}$.

Let us call the icosian $\alpha_{q} q$ the \emph{extension} of the
primitive admissible element $q\in\II$, and $(\alpha^{}_q q,
\alpha^{\,\prime}_{q} \widetilde{q}\,)$ the corresponding
\emph{extension pair}. In view of the form of the rotation $x\mapsto
\frac{1}{\lvert q\widetilde{q}\ts\rvert}\, q x \widetilde{q}$, it is
actually rather natural to replace $q$ and $\widetilde{q}$ by certain
$\oo$-multiples, $q^{}_{\alpha} := \alpha^{}_{q}\ts q$ and
$\widetilde{q^{}_{\alpha}} = \alpha^{}_{\tilde{q}}\,\widetilde{q}$,
such that $\nr(q_{\alpha})$ and $\nr(\widetilde{q_{\alpha}})$ have the
same prime divisors in $\oo$. The introduction of the extension pair
restores some kind of symmetry of the expressions in relation to
the two quaternions involved, which will become evident
in the general treatment of $4$-space \cite{BZ2}.

Clearly, since the extra factors are central, this modification 
does not change the rotation, so that
\begin{equation} \label{unchanged-rot}
     \frac{q x \widetilde{q}}{\lvert q\widetilde{q}\ts\rvert} 
     \, = \, \frac{q_{\alpha} x \widetilde{q_{\alpha}}}
     {\lvert q_{\alpha}\widetilde{q_{\alpha}}\ts\rvert}
\end{equation}
holds for all quaternions $x$.  Note that the definition of the
extension pair is unique up to units of $\oo$, and that one has the
relation
\begin{equation} \label{equal-norms}
     \nr (q_{\alpha}) \, = \,
     \lcm \bigl(\nr (q), \nr(\widetilde{q}\,)\bigr)
     \, = \, \nr (\widetilde{q_{\alpha}})
     \, = \, \lvert q_{\alpha}\ts
             \widetilde{q_{\alpha}}\rvert
     \, \in \, \NN \, ,
\end{equation}
which will be crucial later on.

\begin{lemma} \label{divide} 
  {}For $q\in\II$ and $\gamma\in K$, one has $q\in \gamma \II$
  if and only if\/ $\{ \tr (q\bar{y})\mid y\in\II\}\subset \gamma \oo$. 
\end{lemma}
\begin{proof}
  The statement is clear for $\gamma=0$, so assume $\gamma\neq 0$.
  When $q\in\II$, one has $\tr (q\bar{y})\in\oo$ for all $y\in\II$ (as
  then $q\bar{y}\in\II$), whence $q\in \gamma \II$ implies $\tr
  (q\bar{y})\in \gamma\oo$.  Conversely, $\tr (q\bar{y})\in\oo$ for
  all $y\in\frac{1}{\gamma}\II$ means $q\in
  \bigl(\frac{1}{\gamma}\II\bigr)^* = \gamma \II^* = \gamma \II$, by
  Fact~\ref{self-dual}, which implies the claim.
\end{proof}

\begin{lemma} \label{gcd-trick} 
  If $q\in\II$ is primitive, there is a quaternion $z\in\II$ with $\tr
  (q\bar{z}) = 1$. When, in addition, $q$ is also admissible, there
  exists a quaternion $z\in\II$ such that\/ $\tr
  (q^{}_{\alpha}\ts\bar{z}) + \tr
  (\ts\widetilde{\bar{z}\,}\ts\widetilde{q^{}_{\alpha}}) = 1$, where
  $q^{}_{\alpha}$ denotes the extension of $q$.
\end{lemma}
\begin{proof}
When $q\in\II$, one has $\gcd \{ \tr (q\bar{x})\oo \mid x\in\II\}
=\gamma\oo$ with $\gamma\in\oo$. If $\gamma$ is \emph{not} a unit in
$\oo$, one has $q\in\gamma\II$ by Lemma~\ref{divide}, whence $q$
cannot be primitive in this case. So, $\gamma$ must be a unit,
hence $\gamma\oo=\oo$. Then, by standard arguments based on the
prime ideals, there are \emph{finitely} many icosians
$x_i\in\II$, say $\ell$ of them (in fact, $\ell\le 4$ suffices), such that
\[
    \gcd \{ \tr (q\bar{x}_i)\oo \mid 1\le i \le \ell\}
    \, = \, \tr (q\bar{x}^{}_{1})\oo + \ldots +
    \tr (q\bar{x}^{}_{\ell})\oo \, = \, \oo \, .
\]
This implies the existence of numbers
$\beta_i\in\oo$, with $1\le i\le \ell$, such that $z=\sum_i
\beta_i x_i$ satisfies $\tr (q\bar{z})=1$.

{}For the second claim, assume that $q$ is also admissible and
denote its extension by $q^{}_{\alpha}$. Let $z\in\II$
be the icosian from the first part of the proof, so that
$\tr (q\bar{z})=1$. Since $q^{}_{\alpha}=\alpha^{}_{q}\ts q$ with
$\alpha^{}_{q}\in\oo$, this implies $\tr (q^{}_{\alpha}\bar{z}) =
\alpha^{}_{q}$ and thus also
\[
   \alpha^{\,\prime}_{q} \, = \, \,\widetilde{\!\alpha^{}_{q}\!}\, 
   \, = \, \bigl(\tr (q^{}_{\alpha}\bar{z})\bigr)^{\up} \, = \,
   \tr \bigl(\,\widetilde{\!q^{}_{\alpha}\bar{z}\!}\,\bigr) \, = \,
   \tr \bigl( \ts\widetilde{\bar{z}\,}\ts\widetilde{q^{}_{\alpha}}\bigr).
\]
Since the ideals $\alpha^{}_{q}\ts\oo$ and
$\alpha^{\,\prime}_{q}\ts\oo$ are relatively prime by construction, we
have $\alpha^{}_{q}\ts\oo + \alpha^{\,\prime}_{q}\ts\oo = \oo$ and
thus the existence of $\beta,\delta\in\oo$ with $\beta\alpha^{}_{q} +
\delta\alpha^{\,\prime}_{q} = 1$. The icosians $x=\beta z$ and
$y=\delta{\ts}' z$ then satisfy $\tr (q^{}_{\alpha}\ts \bar{x} ) +
\tr ( \widetilde{\bar{y}\,}\widetilde{q^{}_{\alpha}}\ts) = 1$ as
well as $
\tr ( \widetilde{\bar{x}\,}\widetilde{q^{}_{\alpha}}\ts)
+ \tr (q^{}_{\alpha}\ts \bar{y} ) = 1$, where the second
identity follows from the first via $(\tr (u\bar{v}))^{\up}
= \tr (\ts \widetilde{\bar{v}\,}\widetilde{u}\,)$.

{}Finally, observe that $\tr (u\bar{v})\in K$ for all $u,v\in\HH (K)$,
so that one also has the relation $( \tr (u\bar{v}))^{\prime} = \tr
(\widetilde{\bar{v}\,} \widetilde{u}\ts)$.  Consequently, defining $z =
\tau\ts x + (1-\tau)\ts y$ with the $x,y$ from above, $z$ is an
icosian that satisfies
\[
    \tr (q^{}_{\alpha}\ts\bar{z}) + 
    \tr (\widetilde{\bar{z}\,} \widetilde{q^{}_{\alpha}})
  \, = \, \tau \bigl( \tr (q^{}_{\alpha}\ts\bar{x})
    +  \tr (\widetilde{\bar{y}\,} \widetilde{q^{}_{\alpha}})\bigr)
  + (1-\tau)\bigl(\tr (q^{}_{\alpha}\ts\bar{y})
    +  \tr (\widetilde{\bar{x}\,} \widetilde{q^{}_{\alpha}})\bigr)
  \, = \,  1 \, ,
\]
which establishes the second claim.
\end{proof}

{}For our further discussion, it is convenient to define the set
\begin{equation} \label{sublattices}
     \alat (q) \, = \, \{ qx + \widetilde{x}\widetilde{q}
     \mid x\in\II \} \, = \, \phi^{}_{+} (q\ts\II)\, ,
\end{equation}
which is a sublattice of $\alat$, compare Lemma~\ref{new-sum}.  Note
that, due to $\widetilde{\II}=\II$, one has $\alat (q) =
\widetilde{\alat (q)}$.
\begin{theorem} \label{csl-form}
   Let $q\in\II$ be admissible and primitive,
   and let $q_{\alpha} = \alpha_{q}\ts q$ be its extension.
   Then, the CSL defined by $\ts q$ is given by
\[
   \alat \cap \frac{1}{\lvert q\widetilde{q}\ts\rvert}\ts
    q \alat\widetilde{q} \, = \, \alat (q_{\alpha})\, ,
\]
   with $\alat (q_{\alpha})$ defined as in Eq.~$\eqref{sublattices}$.
\end{theorem}
\begin{proof}
To show the equality claimed, we have to establish two inclusions,
where we may use the fact that $q$ and $q_{\alpha}$ define the
same rotation in $4$-space, see Eq.~\eqref{unchanged-rot}.

First, since $\alat (q_{\alpha})\subset \alat$ is clear, we need to
show that $\lvert q^{}_{\alpha}\widetilde{q^{}_{\alpha}} \rvert\,\alat
(q_{\alpha})\subset q^{}_{\alpha}\alat\widetilde{q^{}_{\alpha}}$. If
$x\in\alat (q_{\alpha})$, there is some $y\in\II$ with
$x=q^{}_{\alpha}\ts y + \widetilde{y} \widetilde{q^{}_{\alpha}}$.
Consequently, observing the norm relations from
Eq.~\eqref{equal-norms}, we find
\[
  \lvert q^{}_{\alpha}\widetilde{q^{}_{\alpha}}\rvert\, x \, = \,
  q^{}_{\alpha}\ts y \,
  \,\overline{\!\widetilde{q^{}_{\alpha}}\!}\,\widetilde{q^{}_{\alpha}}
  + q^{}_{\alpha} \,\overline{\! q^{}_{\alpha} \!}\,
  \widetilde{y}\ts\widetilde{q^{}_{\alpha}} \, = \, q^{}_{\alpha} \,
  (y \,\overline{\!\widetilde{q^{}_{\alpha}}\!}\, + \,\overline{\!
  q^{}_{\alpha}\!}\,\widetilde{y}\ts ) \, \widetilde{q^{}_{\alpha}} \,
  \in \, q^{}_{\alpha} \alat (q^{}_{\alpha})\widetilde{q^{}_{\alpha}}
  \, \subset \, q^{}_{\alpha} \alat \widetilde{q^{}_{\alpha}}\, ,
\]
which gives the first inclusion.

Conversely, when $x\in \alat\cap \frac{1}{\lvert q\widetilde{q}\ts
  \rvert}\ts q \alat\widetilde{q}$, Lemma~\ref{gcd-trick} tells us
  that there exists an icosian $z\in\II$ such that $\tr
  (q_{\alpha}\bar{z}) + \tr (\widetilde{\bar{z}\,}
  \widetilde{q_{\alpha}} ) = 1$. At the same time, there is some
  $y\in\alat$ so that $x=\frac{q y \tilde{q}}{\lvert q
  \tilde{q}\rvert} = \frac{q^{}_{\alpha} y \,\widetilde{\!
  q^{}_{\alpha}\! }} {\lvert q^{}_{\alpha} \,\widetilde{\!
  q^{}_{\alpha}\!} \,\rvert}$.  Observing $x=\widetilde{x}$ and the
  norm relations in Eq.~\eqref{equal-norms}, one finds
\[
   x \, = \, \tr (q^{}_{\alpha}\bar{z})\ts x + \widetilde{x} 
     \tr (\widetilde{\bar{z}\,} \widetilde{q_{\alpha}}) \, = \,
     ( q^{}_{\alpha}\bar{z} + z\,\lbar{q^{}_{\alpha}}) \ts x
    + \widetilde{x}\ts (\widetilde{\bar{z}\,} \widetilde{q_{\alpha}}
    + \lbar{\widetilde{q_{\alpha}}}\widetilde{z}\, ) \, = \,
    q^{}_{\alpha}\ts (\bar{z}x + \widetilde{y}\ts\widetilde{z}\,) +
    (\widetilde{x}\widetilde{\bar{z}\,}\! + zy )
   \ts \widetilde{q^{}_{\alpha}}\, ,
\]
which shows $x$ to be an element of $\alat (q^{}_{\alpha})$.
\end{proof}

\section{Coincidence indices and generating functions}

With the explicit identification of the CSL that emerges from the
rotation defined by an admissible primitive icosian $q$, one can then
calculate the corresponding index. This is either possible by a more
direct (though somewhat tedious) calculation along the lines of
reference \cite{BPR} or by relating $\varSigma^2$ to the corresponding
index of the coincidence site module of $\II$, see \cite{BZ2} for
details. The result reads as follows.

\begin{theorem}  \label{index-formula}
  If $q\in\II$ is an admissible primitive icosian, the rotation
  $x\mapsto \frac{1}{\lvert q\widetilde{q}\ts\rvert}
  qx\widetilde{q}$ is a coincidence isometry of $\alat$.
  Moreover, the corresponding coincidence index satisfies 
\[
   \bigl(\varSigma (q)\bigr)^2 \, = \, \N\bigl(
  \lcm (\nr(q),\nr(q)') \bigr) = \N\bigl(\nr(q^{}_{\alpha})\bigr)\, ,
\] 
  which is a square in $\NN$. Equivalently, one has the formula
\[
   \varSigma(q) \, = \, \nr(q^{}_{\alpha}) \, = \,
   \lcm (\nr(q),\nr(q)') \, ,
\]
which is then, with our above convention from Eq.~$\eqref{def-alpha}$,
always an element of $\ts\NN$.  \qed
\end{theorem}

Due to the subtle aspects of the factorisations of icosians into
irreducible elements, we have not yet found a clear and systematic
approach to the number $f(m)$ of CSLs of $A_4$ of index $m$, though we
will indicate later what the answer might look like. However, at this
point, it is possible to determine the number $f_{\rm rot} (m)$ for
$m$ a prime power.

Clearly, we have $f_{\rm rot}(1)=1$. When $g(m)$ denotes the number of
\emph{primitive} SSLs of $A_4$ of index $m^2$, one can immediately
extract some cases from the explicit results in \cite{H,BHM}. In
particular, one has $f_{\rm rot}(p^r) = g(p^{2r})$ both for $p=5$ and
for all rational primes $p\equiv \pm 2$ mod $5$.  The remaining case
with $p\equiv \pm 1$ mod $5$, where $p$ splits as $p=\pi \pi'$ on the
level of $\oo$, is slightly more difficult, because one has to
keep track of how the algebraically conjugate primes of $\oo$
are distributed between $\nr (q)$ and $\nr(\widetilde{q}\ts)$. Observe
the relation
\[
   \frac{1+p^{-2s}}{1-p^{1-2s}} \, = \,
    1 + \sum_{\ell\ge 1} (p^{\ell} + p^{\ell -1})\ts p^{-2\ell s},
\]
which, for $p\equiv\pm 1$ mod $5$, happens to be the generating
function for the primitive right ideals $q\II$ of the icosian ring of
$p$-power index such that $\nr (q)$ is a power of
$\pi$ (up to units). Those with $\nr (q)$ a power of $\pi'$ produce an
Euler factor of the same form. With this, one can explicitly calculate
$f_{\rm rot}(p^r)$ by collecting all contributions to the index
$\varSigma = p^r$ according to Theorem~\ref{index-formula}. This turns
out to be completely analogous to the calculations for the centred
hypercubic lattice in $4$-space presented\footnote{Note that the
  arithmetic functions in \cite{B,Z1}, in the case of $4$ dimensions, 
  also count the coincidence rotations in multiples of the number of rotation
  symmetries, and \emph{not} the CSLs themselves, hence giving
  the generating function \eqref{def-diri-rot} rather than \eqref{def-diri} in
  this case.}  in \cite{B,Z1}, and the formula for $f$ reads
\begin{equation} \label{mult-fun}
   f_{\rm rot}(p^r) \, = \, \begin{cases}
     6\cdot 5^{2r-1} , & \text{if } p=5  \, , \\
     \frac{p+1}{p-1}\, p^{r-1} (p^{r+1} + p^{r-1} - 2)
      , & \text{if } p\equiv\pm 1\; (5)\, , \\
     p^{2r}+p^{2r-2} , & \text{if } p\equiv\pm 2\; (5)  \, .
   \end{cases}
\end{equation}

\begin{theorem}
  Let $120\ts f_{\rm rot}(m)$ be the number of coincidence rotations of
  $A_4$ of index $m$. Then, $f_{\rm rot}(m)$ is a multiplicative
  arithmetic function, with Dirichlet series generating
  function
\begin{equation} \label{prime-powers}
   \begin{split}
   \varPhi_{A_4}^{\rm \, rot} (s) & \, = \, \sum_{m=1}^{\infty}
   \frac{f_{\rm rot}(m)}{m^s} \, = \, \frac{\zeta^{}_{K} (s-1)}{1+5^{-s}}\,
   \frac{\zeta(s)\,\zeta(s-2)}{\zeta(2s)\,\zeta(2s-2)} \\[2mm]
   & \, = \, \frac{1+5^{1-s}}{1-5^{2-s}}\, \prod_{p\equiv\pm 1 \,(5)}
    \frac{(1+p^{-s})\,(1+p^{1-s})}{(1-p^{1-s})\,(1-p^{2-s})}\,
    \prod_{p\equiv\pm 2 \,(5)} \frac{1+p^{-s}}{1-p^{2-s}}  \\[2mm]
    & \, = \, 1 + \frac{5}{2^s} + \frac{10}{3^s} + \frac{20}{4^s} + 
      \frac{30}{5^s} + \frac{50}{6^s} + \frac{50}{7^s} + \frac{80}{8^s}
        + \frac{90}{9^s} + \frac{150}{10^s} + \frac{144}{11^s} + \ldots
   \end{split}
\end{equation}
  where $\zeta(s)$ is Riemann's zeta function and $\zeta^{}_{K} (s)$
  denotes the Dedekind zeta function of the quadratic field 
  $K=\QQ(\sqrt{5}\,)$.
\end{theorem}
\begin{proof}
  The multiplicativity of $f_{\rm rot}$ is inherited from the unique
  factorisation in $\II$ together with the divisor properties of the
  coincidence index from Lemma~\ref{factor-lemma}.  Consequently,
  Eq.~\eqref{prime-powers} fixes $f_{\rm rot}(m)$ for all $m\in\NN$
  via the identity $f_{\rm rot}(mn)=f_{\rm rot}(m)\ts f_{\rm rot}(n)$
  for integers $m$ and $n$ that are relatively prime.

It is a routine exercise to calculate the Euler factors of the
corresponding Dirichlet series generating function and to express
$\varPhi_{A_4}^{\rm \, rot}$ in terms of the two zeta functions mentioned.
\end{proof}

Note that $\zeta^{}_{K} (s) = \zeta(s)\, L(s,\chi)$, where $L(s,\chi)$
is the $L$-series of the primitive Dirichlet character $\chi$ defined
by
\[
     \chi(n) \, = \, \begin{cases}
      0, &  n\equiv 0 \; (5)\, , \\
      1, &  n\equiv \pm 1 \; (5)\, , \\
     -1, &  n\equiv \pm 2 \; (5)\, . \end{cases}
\]

Observe next that $f_{\rm rot}(m)>0$ for all $m\in\NN$, so that also
the number $f(m)$ of CSLs must be positive (though we can still have
$0<f(m)<f_{\rm rot} (m)$). Moreover, each element of $\OC (A_4)$ can
be written as a product of a rotation with a reflection that maps
$A_4$ onto itself. Consequently, the simple coincidence spectrum
$\varSigma (\OC (A_4))$ is the set of all positive integers. Although
multiple coincidences may produce further lattices, compare
\cite{BG,Z2,BZ1}, the total spectrum $\varSigma_{A_4}$ cannot be
larger than the elementary one, so that the following consequence is
clear.

\begin{coro} \label{spectra}
   The multiple coincidence spectra of the root lattice $A_4$
   coincide with the elementary one, and one has $\varSigma_{A_4} =
   \varSigma (\OC (A_4)) = \varSigma (\SOC (A_4)) = \NN$, which is a
   monoid. \qed
\end{coro}

The Dirichlet series $ \varPhi_{A_4}^{\rm \, rot} (s)$ is analytic in the
open right half-plane $\{s=\sigma+ \ii t \mid \sigma > 3\}$, and has a
simple pole at $s=3$. The corresponding residue is given by
\begin{equation} \label{residue}
   \mathrm{res}_{s=3}\, \varPhi_{A_4}^{\rm \, rot} (s) \, = \,
   \frac{125}{126}\,\frac{\zeta^{}_{K} (2)\,\zeta(3)}{\zeta(4)\,\zeta(6)}
   \, = \, \frac{450 \sqrt{5}}{\pi^6}\, \zeta(3) \, \simeq \,
   1.258{\,}124\, ,
\end{equation}
which is based on the special values
\[
    \zeta(4) \, = \, \frac{\pi^4}{90} \; , \quad
    \zeta(6) \, = \, \frac{\pi^6}{945} \; , \quad
    \zeta^{}_{K} (2) \, = \, \frac{2\ts\pi^4}{75\sqrt{5}}\; ,
\]
together with $\zeta (3) \simeq 1.202{\,}057$, compare \cite{BM}
and references given there. The value of $\zeta(3)$ is known to be
irrational, but has to be calculated numerically.

With this information, we can extract the asymptotic behaviour of the
counts $f_{\rm rot}(m)$ from the generating function
$\varPhi_{A_4}^{\rm \, rot} (s)$ by Delange's theorem, see
\cite[Appendix]{BM} for a formulation tailored to this situation.  One
obtains, as $x\to\infty$,
\begin{equation} \label{asymp}
   \sum_{m\le x} f_{\rm rot}(m) \, \sim \, 
   \mathrm{res}_{s=3}\, \varPhi_{A_4}^{\rm \, rot} (s)
   \, \frac{x^3}{3} \, \simeq \, 0.419{\,}375\, x^3 .
\end{equation}
Clearly, this is also an upper bound for the asymptotic behaviour
of the true CSL counts. 

As mentioned earlier, one is primarily interested in the number $f(m)$
of CSLs of index $m$, which satisfies $0 < f(m) \le f_{\rm rot} (m)$ for all
$m\in\NN$ in view of Corollary~\ref{spectra}. Some preliminary 
calculations show that also $f(m)$ is multiplicative. In fact, one has
$f(p^r) = f_{\rm rot}(p^r)$ for all primes $p \equiv \pm 2 \; (5)$, and
$f(5^r) = f_{\rm rot}(5^r)/5$ for $r\ge 1$. For the remaining primes
$p\equiv\pm 1\; (5)$, one has $f(p)=f_{\rm rot}(p)$, but differences 
occur for all powers $p^r$ with $r\ge 2$. This happens first for
$m=11^2=121$ and is induced by the more complicated factorisation
for these primes, compare \cite{BG-shell} for a similar phenomenon.
Consequently, the modification for the prime $5$ is 
sufficient up to index $m=120$. The Dirichlet
series generating function thus starts as
\[
 \varPhi_{A_4}^{} (s) \, = \, \sum_{m=1}^{\infty}
   \frac{f(m)}{m^s} 
   \, = \,   1 + \frac{5}{2^s} + \frac{10}{3^s} + \frac{20}{4^s} + 
      \frac{6}{5^s} + \frac{50}{6^s} + \frac{50}{7^s} + \frac{80}{8^s}
        + \frac{90}{9^s} + \frac{30}{10^s} + \frac{144}{11^s} + \ldots     
\]
At this stage, the general mechanism behind this is not completely
unravelled, but we hope to present it in greater generality in \cite{BZ2}.

\section{Related results and outlook}

In one dimension, the CSL problem becomes trivial, so that
$\varPhi(s)\equiv 1$ in this case. In two dimensions, a rather general
approach to lattices and modules is possible via classic algebraic
number theory, see \cite{PBR,BG} and references therein, which
includes the treatment of multiple coincidences. For the root lattice
$A_2$, the CSL generating function reads
\begin{equation} \label{planar-case}
  \varPhi_{A_2} (s) \, = \, \varPhi_{A_2}^{\rm \, rot} (s) 
     \, = \, \prod_{p\equiv 1 \; (3)}
     \frac{1+p^{-s}}{1-p^{-s}} \, = \, \frac{1}{1+3^{-s}}\,
     \frac{\zeta^{}_{\QQ(\xi_3)}(s)}{\zeta(2s)}\, ,
\end{equation}
where $\xi_3$ is a primitive cube root of $1$. Here, the equality is a
consequence of the commutativity of $\SOC (A_2)$.  The simple
coincidence spectrum of this lattice is the multiplicative monoid of
integers that is generated by the rational primes $p\equiv 1$ mod $3$.

In $3$-space, various examples are derived in \cite{B} and have
recently been proved by quaternionic methods \cite{BPR} similar
to the ones used here. Among these cases is the root lattice $A_3$,
which happens to be the face centred cubic lattice in $3$-space,
with generating function
\begin{equation} \label{three-space}
   \varPhi_{A_3} (s) \, = \,  \varPhi_{A_3}^{\rm \, rot} (s) 
     \, = \,  \prod_{p\neq 2}
     \frac{1+p^{-s}}{1-p^{1-s}} \, = \,  \frac{1-2^{1-s}}{1+2^{-s}}\,
     \frac{\zeta(s)\,\zeta(s-1)}{\zeta(2s)}\, .
\end{equation}
The equality of the two Dirichlet series to the left is non-trivial,
and was proved in \cite{BPR} with an argument involving Eichler
orders. The same formula also applies to the other cubic lattices in
$3$-space \cite{B,BZ1}.  The simple coincidence spectrum is thus the
set of odd integers, which is again a monoid. The multiple analogues
have recently been derived in \cite{Z2,Z3}, see also \cite{BG,BZ1} for
related results.

Several of these results are also included by now in \cite{online}.
In $4$-space, various other lattices and modules of interest exist,
for which some results are given in \cite{B,Z1,Z3}, with more structural
proofs and generalisations being in preparation \cite{BZ2}. Beyond
dimension $4$, very little is known \cite{Zou1,Zou2}, though it should
be possible to derive the simple coincidence spectra for certain
classes of lattices.

\smallskip
\section*{Acknowledgements}
It is our pleasure to thank Johannes Roth for his cooperation, and
Robert V.\ Moody and Rudolf Scharlau for helpful discussions.  This
work was supported by the German Research Council (DFG), within the
CRC 701, and by EPSRC, via grant EP/D058465.

\end{document}